\documentclass[journal,twoside,web]{ieeecolor}
\usepackage{lcsys}
\usepackage{amsmath,amsfonts,mathtools,amssymb}
\usepackage{graphicx}
\usepackage{textcomp}
\usepackage{tikz}
\usetikzlibrary{arrows.meta,calc,decorations.markings,math}
\usepackage{pgfplots}
\usepackage{aligned-overset}
\usepackage{nicefrac}
\usepackage{algorithm}
\usepackage{lipsum}
\usepackage{soul}
\usepackage{siunitx}
\DeclareSIUnit{\byte}{B}
\usetikzlibrary{positioning}
\usepackage[hidelinks]{hyperref}
\allowdisplaybreaks
\usepackage[noend]{algpseudocode}
\usepackage{yfonts}
\usepackage{tabularx}
\usepackage{cite}
\let\labelindent\relax
\usepackage[shortlabels]{enumitem}
\usepackage{comment}
\usepackage{changepage} % in preamble
\usepackage{lipsum}
\usepackage[normalem]{ulem}

\newcommand{\safeopt}{\textsc{SafeOpt}}

\newcommand{\ie}{i\/.\/e\/.,\/~}
\newcommand{\eg}{e\/.\/g\/.,\/~}

\newtheorem{assumption}{\bf Assumption}
\newcommand{\domain}{\mathcal A}
\newcommand{\E}{\mathbb{E}}
\renewcommand{\P}{\mathbb{P}}

\newcommand{\AT}[1]{\textcolor{black}{#1}}
\newcommand{\fakepar}[1]{\vspace{1mm}\noindent\textbf{#1.}}
\DeclareMathOperator*{\argmax}{arg\,max\,}

\definecolor{aaltoRed}{RGB}{239,51,64}%
\definecolor{aaltoBlue}{RGB}{0,94,184}%

\DeclareMathOperator*{\R}{\mathbb{R}}

\newcounter{tool}

\definecolor{magenta}{RGB}{255,0,255}
\definecolor{lightgreen}{RGB}{76, 200, 76}
\definecolor{deeporange}{RGB}{204, 112, 0}

\newcommand{\I}{\mathcal{I}}
\newcommand{\Ig}{\mathcal{I}_{\mathrm{g}}}

\newtheorem{lemma}{\bf Lemma}
\newtheorem{theorem}{\bf Theorem}

\newtheorem{corollary}{\bf Corollary}

\newtheorem{remark}{\bf Remark}

\title{Safe Bayesian optimization across noise models via scenario programming}

\author{Abdullah Tokmak, Thomas B.\ Schön, and Dominik Baumann
\thanks{This research was partially supported by \emph{Kjell och Märta Beijer Foundation.}}
\thanks{Abdullah Tokmak and Dominik Baumann are with the Cyber-physical systems group, Aalto University, Espoo, Finland (email: \textit{firstname.lastname@aalto.fi}).}
\thanks{Thomas B.\ Schön is with the Department of Information Technology, Uppsala University, Uppsala, Sweden (email: \textit{thomas.schon@uu.se}).}
}

\usepackage{fancyhdr}
\fancypagestyle{firstpage}{%
  \fancyhf{} % clear header/footer
  \fancyfoot[L]{\normalfont \sffamily \scriptsize \lowernote}
}
\fancypagestyle{plain}{%
  \fancyhf{}
  \fancyfoot[L]{\normalfont \sffamily \scriptsize \lowernote}
}

\newcommand{\lowernote}{
\textbf{Accepted final version. \textcopyright\; 2025
IEEE.}
To be published at IEEE Control Systems Letters.
Personal use of this material is permitted. Permission from IEEE must be obtained for all other uses, in any current or future media,
including reprinting/republishing this material for advertising or promotional purposes, creating new collective works, for resale or redistribution to servers
or lists, or reuse of any copyrighted component of this work in other works.
}
\begin{document}
\maketitle
\thispagestyle{firstpage}
% configurable
\begin{abstract}
    Safe Bayesian optimization (BO) with Gaussian processes is an effective tool for tuning control policies in safety-critical real-world systems, specifically due to its sample efficiency and safety guarantees.
    However, most safe BO algorithms assume homoscedastic sub-Gaussian measurement noise, an assumption that does not hold in many relevant applications.
    In this article, we propose a straightforward yet rigorous approach for safe BO across noise models, including homoscedastic sub-Gaussian and heteroscedastic heavy-tailed distributions.
    We provide a high-probability bound on  the measurement noise  via the scenario approach, integrate these bounds into high probability confidence intervals, and prove safety and optimality for our proposed safe BO algorithm.
We deploy our algorithm in synthetic examples and in tuning a controller for the Franka Emika manipulator in simulation.
\end{abstract}

\section{Introduction}
% Tuning parametrized control policies is a constructive approach to improve the performance of many robotic systems.
Learning-based methods, such as reinforcement learning (RL)~\cite{sutton2018reinforcement}, have been successfully used to obtain high-performing policies for nonlinear systems with unknown dynamics.
However, to tune control policies of safety-critical real-world systems, two essential criteria must be met.
First, we require sample efficiency, since obtaining a sample corresponds to conducting a real-world experiment.
Second, we need safety guarantees to ensure secure operation of the system.
Vanilla RL has inherent difficulties with both criteria.
Bayesian optimization (BO)~\cite{garnett2023bayesian} with Gaussian process (GP)~\cite{Rasmussen2006Gaussian} surrogates offers a sample-efficient alternative that works effectively in, \eg robotics applications~\cite{calandra2016bayesian}.
Moreover, safe BO algorithms~\cite{sui2015safe} also provide probabilistic safety guarantees. 
% However, most safe BO works assume that the observation noise is homoscedastic sub-Gaussian.
However, most safe BO works assume---and are thus only robust with respect to---homoscedastic~$R$-sub-Gaussian measurement noise.
This assumption is broken in many relevant scenarios, such as in the modeling of network delays~\cite{jagannathan2014throughput},  measurements from radar or LiDAR sensors~\cite{zou2024vb}, or parameter tuning on real-world systems, where the observation noise often exhibits strong dependence on the chosen parameter~\cite{kirschner2018information, makarova2021risk}.
In this letter, we overcome this bottleneck by proposing a practical, safe BO algorithm that is straightforwardly configurable across different noise models. 

\fakepar{Related work}
\safeopt~\cite{sui2015safe}---which extends GP-UCB~\cite{srinivas2012information, chowdhury2017kernelized} to probabilistic constraint satisfaction---is arguably the most popular safe BO algorithm.
\safeopt\ and extensions have been applied to control parameter tuning for medical devices~\cite{sui2015safe}, quadrotor vehicles~\cite{berkenkamp2023bayesian}, the Franka Emika manipulator~\cite{Bhavi2023GSO}, and the Furuta pendulum~\cite{tokmak2025safe}.
Although the seminal \safeopt\ work~\cite{sui2015safe} requires independent and identically distributed (i.i.d.) Gaussian noise~$\epsilon_t$,~\cite{tokmak2024pacsbo, tokmak2025safe, berkenkamp2023bayesian, Bhavi2023GSO} relax this assumption to~$\epsilon_t$ only being homoscedastic $R$-sub-Gaussian conditioned on the filtration, \ie for all~$\lambda \in \R,$ $\E[\exp(\lambda \epsilon_t)\vert \mathcal F_{t-1}]\leq \exp(\nicefrac{\lambda^2 R^2}{2})$ for all iterations~$t\geq 1$.
In the sub-Gaussian case, the safe BO algorithms use confidence intervals derived in, \eg~\cite{chowdhury2017kernelized, 
Fiedler2021Practical}.

Nevertheless, many relevant applications~\cite{zou2024vb, jagannathan2014throughput} are endowed by heavy-tailed distributions whose tails do not decay sufficiently fast to be classified as sub-Gaussian.
In the linear bandit setting,~\cite{bubeck2013bandits,shao2018almost} propose algorithms capable of handling heavy-tailed noise, while~\cite{tajdini2025improved} extends these approaches to the kernelized setting.
These works only require a (known) bound on the $(1+\alpha)$-th moment for~$\alpha\in(0,1]$.
Moreover,~\cite{chowdhoury2019bayesian} proposes a BO algorithm based on GP-UCB for heavy-tailed rewards.
%Specifically,~\cite{chowdhoury2019bayesian} adaptively truncates the reward into a finite-dimensional approximate feature space by using quadrature Fourier features or Nyström approximations.
However, the mentioned works~\cite{bubeck2013bandits, shao2018almost, tajdini2025improved, chowdhoury2019bayesian} only consider homoscedastic noise and do not account for constraints. % satisfaction.

References~\cite{kirschner2018information, makarova2021risk} argue that many use cases involve  heteroscedastic noise, \ie noise that depends on the input.
To address such cases,~\cite{kirschner2018information} models the stochastic bandit problem under heteroscedastic sub-Gaussian noise and derives corresponding confidence intervals for functions in reproducing kernel Hilbert spaces (RKHSs).
This work is then extended in~\cite{makarova2021risk} to the risk-averse BO setting.

While the contributions mentioned above tackle specific noise assumptions, none of these works are applied to safe BO or account for heteroscedastic heavy-tailed distributions.
More crucially, each noise assumption relies on its own theoretical foundation and demands intricate transformations for implementation.
Altogether, there is no unified and practical yet rigorous framework for safe BO that is straightforwardly adjustable across noise models, including homoscedastic sub-Gaussian and heteroscedastic heavy-tailed distributions.

\fakepar{Contribution}
We propose a safe BO algorithm reminiscent of \safeopt\  that is straightforwardly configurable across different noise models.
We make the following contributions:\footnote{%
    Code/experiments: \url{https://github.com/tokmaka1/ACC-2026}}
\begin{adjustwidth}{0.5cm}{}
\begin{enumerate}[label=(C\arabic*)]
    \item We propose a scenario-based approach~\cite{campi2018introduction} to obtain probabilistic noise bounds.  %bound the noise using sampling-based uncertainty quantification via the scenario approach~\cite{campi2018introduction} and  \AT{derive} high-probability confidence intervals. 
    \label{con:noise}
    \item Using the bounds from~\ref{con:noise}, we
    derive high probability confidence intervals. \label{con:confidence}
    \item  We develop a safe BO algorithm reminiscent of \safeopt~\cite{sui2015safe} using the confidence intervals from~\ref{con:confidence}, for which we prove safety and optimality. \label{con:safety_optimality}
\end{enumerate} 
\end{adjustwidth}

\section{Problem setting and preliminaries}
Next, we explain the optimization problem (Sec.~\ref{sec:opt}), \safeopt\ (Sec.~\ref{sec:safeopt}), and the scenario approach (Sec.~\ref{sec:scenario}).

\subsection{Optimization problem}\label{sec:opt}
% We work with parametrized policies.
% Hence, we can cast the policy optimization problem as a constrained optimization problem.
We seek to maximize the unknown reward function $f:\mathcal A\subseteq\R^n\rightarrow \R$ subject to unknown constraints $g_i:\mathcal A\subseteq\R^n\rightarrow \R, i\in\Ig\subseteq\mathbb N$.
Analogous to~\cite{berkenkamp2023bayesian}, we introduce a surrogate function~$h_i$ with~$h_i(\cdot)=f(\cdot)$ if~$i=0$ and~$h_i(\cdot)=g_i(\cdot)$ if~$i\in\Ig$, and define~$\I\coloneqq \{0\}\cup\Ig$.
We assume that at each iteration~$t\geq 1$, an experiment following a parametrized policy with policy parameter~$a_t\in\domain$ yields a noisy evaluation $y_{i,t} = h_i(a_t)+\epsilon_{i,t}$, and denote the observation noise vector by $\epsilon_t\coloneqq [\epsilon_{0,t},\ldots,\epsilon_{\lvert\Ig\rvert,t}]^\top$.
\begin{assumption}\label{asm:noise_space}
    The observation noise~$\epsilon_t$ is defined on the probability space $(\Omega, \mathcal F, \mathbb P)$, from which we can sample.
\end{assumption}
% Assumption~\ref{asm:noise_space}  introduces a general and unified view on the noise and allows us to model \emph{any} noise distribution as long as we can sample from its probability space.
Assumption~\ref{asm:noise_space}  introduces a fairly general and unified approach for any noise model, as long as we can sample from its probability space. %; we elaborate on this in Remark~\ref{re:noise}.}
In contrast to the classic~$R$-sub-Gaussian assumption~\cite{chowdhury2017kernelized}, we assume that the noise can be sampled directly, which is not possible under the sub-Gaussian formulation since it represents a family of distributions rather than a specific one.
This assumption enables computing high-probability noise bounds, as established later in
Theorem~\ref{th:simultaneous}.
Arguably, the tail behavior of the noise is the most important aspect for achieving such bounds.
Since the normal distribution~$\mathcal N(0,R^2)$ saturates the moment-generating-function of the $R$-sub-Gaussian family and therefore exhibits the heaviest admissible tails, it provides an empirical and practically convenient conservative surrogate. % for implementing the most conservative case.
Assumption~\ref{asm:noise_space} thereby unifies different noise models under a single theory, allowing practitioners to incorporate prior knowledge about the noise or to employ data-driven or oracle-based representations of the aleatoric uncertainty.

%Altogether,
Altogether, the optimization problem is described by
\begin{align}\label{eq:opt_problem}
    \max_{a\in\domain} f(a)  \;\mathrm{s.t.}\; g_i(a_t) \geq 0, \quad \forall t\geq 1, \forall i \in \Ig,
\end{align}
\ie we sample~$h$ episodically to maximize~$f$ while requiring safety.
In~\eqref{eq:opt_problem}, safety is defined as only conducting experiments that correspond to positive constraint values.
We solve~\eqref{eq:opt_problem} using safe BO to exploit its sample efficiency and safety guarantees by making a smoothness assumption, which is standard in many BO algorithms, see, \eg\cite{sui2015safe,srinivas2012information,chowdhury2017kernelized, berkenkamp2023bayesian, Bhavi2023GSO}.
\begin{assumption}\label{asm:f}
For any~$i\in\I$, $h_i$ is a member of the RKHS of the continuous kernel~$k$ with known  RKHS norm~$\|h_i\|_k$.
\end{assumption}
Note that, for brevity, we assume $k\equiv k_i, \forall i\in I$.
We introduce the kernel metric $d_k(a,a^\prime)\coloneqq \sqrt{2(k(a,a)-k(a,a^\prime))}$ to leverage a well-known result on the continuity of RKHS functions, which will assist us in formulating our safe BO algorithm. % to solve~\eqref{eq:opt_problem}.
\begin{lemma}[{Continuity~\cite[Lemma~1]{Tokmak2023Arxiv}}]\label{le:cont}
Under Assumption~\ref{asm:f},
$
    \lvert h_i(a)-h_i(a^\prime)\rvert \leq \|h_i\|_k  d_k(a,a^\prime),
$
~$\forall a,a^\prime\in\domain, i \in \I$.
\end{lemma}

\subsection{\safeopt}\label{sec:safeopt}
Akin to \safeopt~\cite{berkenkamp2023bayesian}, we use independent GP surrogates with mean function~$\mu_{t,i}$ and variance~$\sigma_{t}$ constructed using kernel~$k$ to model the unknown functions~$h_i$.
To guarantee safety while maximizing the reward, \ie to solve~\eqref{eq:opt_problem}, we must quantify the uncertainty around our belief of the ground truth~$h_i$.
Hence, we require high probability frequentist confidence intervals, %~$h_{i}(a)\in Q_{i,t}(a)$
% of the form 
% \begin{align}\label{eq:Q}
%     Q_{t,i}(a)\coloneqq [\mu_{t,i}(a)\pm \beta_{t,i}\sigma_{t,i}(a)]. 
% \end{align}
% \begin{lemma}[Sub-Gaussian confidence intervals]\label{le:confidence_subgaussian}
% Let Assumption~\ref{asm:f} and~$\epsilon_{i,t}$ be conditionally $R$-sub-Gaussian.
% Then, jointly for all~$t\geq 1, a\in\domain, i\in\I$, with probability~$1-\delta$,
% \begin{align*}
%     &\lvert h_i(a)-\mu_{t,i}(a)\rvert \leq \beta_{i,t}\sigma_{i,t}(a) \\
%     &\beta_{i,t}\coloneqq \lVert h_i\rVert_k +\frac{R}{\sqrt \eta}\sqrt{
% \log\det\left(\frac{1}{\eta}K_t+I_t\right)-
% \frac{2}{\lvert \mathcal I \rvert}\log(\delta)
%     }.
% \end{align*}
% \end{lemma}
% \begin{proof}
%     The approximation error is divided into a bias term and a noise term via the triangle inequality.
% To bound the bias term, see \eg~\cite[Theorem~2]{chowdhury2017kernelized} or~\cite[Remark~3.13]{abbasi2013online}.
% The noise term is upper-bounded in~\cite[Theorem~1]{chowdhury2017kernelized}, while~\cite[Theorem~1]{Fiedler2021Practical} adapts the bound to be purely data-dependent and to hold for regularization parameters~$\eta\leq 1$.
% We allow for multiple constraints by scaling~$\delta$ with~$\lvert \mathcal I\rvert$.
% \end{proof}
as derived in~\cite{chowdhury2017kernelized, Fiedler2021Practical} for homoscedastic sub-Gaussian measurement noise,
of the form
%Specifically, the confidence intervals for all~$t\geq 1, i\in\I, a\in\domain$ are of the form
\begin{align}\label{eq:C}
C_{i,t}(a) \coloneqq C_{i,t-1}(a) \cap [\mu_{i,t}(a)\pm\beta_{i,t}\sigma_{t}(a)],
\end{align}
for all~$t\geq 1, i\in\I, a\in\domain$, with $C_{i,0}(a)\coloneqq (-\infty,\infty)$ and a safety parameter~$\beta_{i,t}$, which we will later make precise as we derive our bounds.
Specifically, we will present confidence intervals that account for general noise models under Assumption~\ref{asm:noise_space}.
We introduce the upper and lower confidence bounds as $u_{i,t}(a)\coloneqq \max C_t(a)$ and $\ell_{t,i}(a)\coloneqq \min C_t(a)$, respectively, denote the uncertainty width by~$w_{i,t}(a)\coloneqq u_{i,t}(a)-\ell_{i,t}(a),$ and define a safe set
\begin{align}\label{eq:safe_set}
    S_t\coloneqq \cap_{i\in\Ig}\cup_{a\in S_{t-1}} \{
    a^\prime \in\domain \vert \ell_ {t,i}(a)-\|h_i\|_k d_k(a,a^\prime)\geq 0
    \}.
\end{align}
The safe set only contains parameters that satisfy the constraints in~\eqref{eq:opt_problem} with high probability~\cite[Theorem~1]{sui2015safe}.
To start exploration, we assume that $\emptyset\neq S_0 \subseteq\domain$, \ie a non-empty set of initial safe samples is given.
To efficiently balance exploration and exploitation while guaranteeing safety, we introduce the set of potential maximizers
$M_t \coloneqq \{a\in S_t\vert u_{t,0} \geq \max_{a^\prime\in S_t} \ell_{t,0}(a^\prime)\}$
and the set of potential expanders 
   $G_t \coloneqq \{a\in S_t\vert  e_t(a)> 0\} , 
    e_t(a)=\lvert\{b \in \domain\setminus S_t \vert \exists i \in \Ig, u_{t,i}(a)-\|g_i\|_k d_k(a,b)\geq 0\}\rvert.$
The set~$M_t$ contains samples that we believe to safely maximize the reward, while an evaluation of the samples within~$G_t$ may safely expand the safe set. %~\eqref{eq:safe_set}.
% \begin{align}\label{eq:maximizer}
%   M_t \coloneqq \{a\in S_t\vert u_{t,0} \geq \max_{a^\prime\in S_t} \ell_{t,0}(a^\prime)\}  
% \end{align}
% and the set of potential expanders
% \begin{align}\label{eq:expander}
%    G_t &\coloneqq \{a\in S_t\vert  e_t(a)> 0\}  \\
%     e_t(a)&=\lvert\{a^\prime \in \domain\setminus S_t \vert \exists i \in \Ig: u_{t,i}(a)-\|g_i\|_k d_k(a,a^\prime)\geq 0\}\rvert \nonumber.
% \end{align}
Akin to \safeopt, the acquisition function is $a_{t+1}=\argmax_{a_i\in (M_{t}\cup G_{t})}\max_{i\in\I} w_{t,i}(a_i)$, \ie we conduct the next experiment with the most uncertain parameter in~$M_t\cup G_t$.

\subsection{Scenario approach}\label{sec:scenario}
Throughout this article, we derive confidence intervals~$C_{i,t}$ (see~\eqref{eq:C}) that are adjustable across different noise models by leveraging Assumption~\ref{asm:noise_space} and the scenario approach~\cite{campi2018introduction}. %,

Let us denote by~$\tilde\epsilon_t^{(j)}, j\in [1,m_t]$, the \emph{scenarios}, which are drawn i.i.d.\ from the probability space $(\Omega, \mathcal F, \P)$. 
Using the scenarios~$\tilde\epsilon_t^{(j)}$, we aim to obtain a high probability scenario bound~$\bar\epsilon_{t}\coloneqq [\bar\epsilon_{0,t},\ldots,\bar\epsilon_{\lvert\Ig\rvert,t}]^\top$ on the unknown observation noise~$\epsilon_t$ by solving the convex scenario program
\begin{align}\label{eq:opt_scenario}
     \min_{\bar\epsilon_t\in\R_{\geq 0}^{\lvert \I\rvert}} \mathbf{1}^\top \bar \epsilon_t \quad \mathrm{s.t.}\quad \bar\epsilon_{i,t}\geq \lvert \tilde\epsilon_{i,t}^{(j)}\rvert,\quad\forall j\in[1,m_t],
\end{align}
where~$\mathbf 1\coloneqq [1,\ldots,1]^\top$. Since we minimize~$\bar\epsilon_t$ 
 element-wise,
\begin{align}\label{eq:opt_scenario_solution}
    \bar\epsilon_{t} = \big[\max\nolimits_{j\in[1,m_t]}\lvert\tilde \epsilon_{0,t}^{(j)}\rvert, \ldots, \max\nolimits_{j\in[1,m_t]}\lvert\tilde \epsilon_{\lvert \Ig\rvert,t}^{(j)}\rvert\big]^\top
\end{align}
is the solution to~\eqref{eq:opt_scenario}.
Remarkably, by solving~\eqref{eq:opt_scenario}, \ie by generating a finite number of scenarios, the scenario approach generalizes to the true random variable~$\epsilon_t$ as long as it lives on the same probability space.
To quantify the generalization property, we define the violation probability as under-estimating the absolute value of the observation noise, \ie
\begin{align}\label{eq:V}
    V_t(\bar\epsilon_t)\coloneqq \P[\exists i\in\I: \bar\epsilon_{i,t} < \lvert \epsilon_{i,t}\rvert].
\end{align}
The following lemma is a classic result from scenario theory and bounds the violation probability~$V_t(\bar\epsilon_t)$.

\begin{lemma}[{Scenario approach~\cite[Theorem~3.7]{campi2018introduction}}]\label{th:noise}
Fix any~$t\geq 1$, 
let Assumption~\ref{asm:noise_space} hold, and~$\tilde\epsilon_t^{(j)}$, $j\in[1,m_t]$ be i.i.d.\ scenarios from $(\Omega,\mathcal F,\P)$.
Then, for any~$\nu,\kappa\in(0,1)$, if~$m_t$ is such that $
\sum_{s=0}^{\lvert\I\rvert-1}{m_t\choose s}\nu^s(1-\nu)^{m_t-s}\leq \kappa$,
% \begin{align}\label{eq:lemma_2}
%     \P^{m_t}[V_t(\bar\epsilon_t)> \nu]\leq \sum_{s=0}^{\lvert\I\rvert-1}{m_t\choose s}\nu^s(1-\nu)^{m_t-s},
% \end{align}
then $\P^{m_t}[V_t(\bar\epsilon_t)> \nu]\leq \kappa$,
with~$\bar\epsilon_t$ as in~\eqref{eq:opt_scenario_solution}. % and~$\P^{m_t}$ is the product probability space of~$\P$.
\end{lemma}

Lemma~\ref{th:noise} states for any fixed iteration~$t\geq 1$ that the violation probability~\eqref{eq:V} is larger than~$\nu$ with confidence of at most~$\kappa$.
Specifically, it quantifies the probability in the product probability space $(\Omega^{m_t}, \mathcal F^{m_t}, \P^{m_t})$,
%$\mathcal F^m\coloneqq \otimes_{i=1}^{m_t}\mathcal F$,
which naturally arises from the~$m_t$ scenarios being i.i.d.\ samples from $(\Omega,\mathcal F, \P)$. %~\cite[Section~II.A]{romao2022exact}.

% For simplicity of exposition, we only consider the basic scenario approach without discarded constraints, as \eg stated in~\cite{romao2022exact}.
% Using the sampling-and-discarding scenario approach follows mutatis mutandis and can yield slightly tighter bounds.
% \AT{Consider removing for space reasons}

\section{Framework for bounding general noise}
% \AT{In this section, we derive probabilistic bounds that are adjustable across different noise models to establish the confidence intervals~\eqref{eq:C}.}
We derive probabilistic bounds on the measurement noise that hold simultaneously for all iterations using the scenario approach (Sec.~\ref{sec:scenario_uniform}) before constructing confidence intervals that are adjustable across different noise models (Sec.~\ref{sec:confidence}).

\subsection{Scenario-based simultaneous bounds}\label{sec:scenario_uniform}

Although Lemma~\ref{th:noise} provides a scenario-based high-probability noise bound, it only captures a \emph{fixed} time step~$t$. 
To provide guarantees on the measurement noise while solving~\eqref{eq:opt_problem}, we require probabilistic bounds that hold \emph{simultaneously} for all~$t\geq 1$.
To this end, we present Algorithm~\ref{alg:noise}.

\begin{algorithm}
\begin{algorithmic}[1]
    \Require $\kappa$, $\nu$, $(\Omega, \mathcal F, \mathbb P)$, $t$, $\lvert\I\rvert$
    \State $\kappa_t\gets\nicefrac{6\kappa}{\pi^2 t^2}$ 
    \State $m_t \gets \min_{m_t \in\mathbb N}\;\mathrm{s.t.}\; \sum_{s=0}^{\lvert \I\rvert -1}{m_t\choose s}\nu^s(1-\nu)^{m_t-s}\leq \kappa_t$
    \State Generate~$m_t$ i.i.d.\ scenarios~$\epsilon_j$ from $(\Omega, \mathcal F, \mathbb P)$
    \State \Return $\bar\epsilon_t\coloneqq [\bar\epsilon_{0,t},\ldots,\bar\epsilon_{\lvert 
\Ig\rvert,t}]^\top$ \Comment{\eqref{eq:opt_scenario_solution}}
\end{algorithmic}
    \caption{Bounding general noise}
    \label{alg:noise}
\end{algorithm}
Specifically, we want the violation probability~\eqref{eq:V} to be bounded by a user-chosen confidence level~$\kappa$ simultaneously for all iterations~$t$.
Algorithm~\ref{alg:noise} computes the number of scenarios~$m_t$ ($\ell.$~2) based on an iteration-adjusted confidence level~$\kappa_t$ ($\ell.$~1).
Then, it generates the scenarios ($\ell.$~3) and returns the bound~$\bar\epsilon_t$~\eqref{eq:opt_scenario_solution},  which solves the scenario program~\eqref{eq:opt_scenario}.
In Theorem~\ref{th:simultaneous}, we will prove that satisfying the iteration-adjusted confidence level~$\kappa_t$ at each iteration indeed provides a scenario bound~$\bar\epsilon_t$ that bounds~\eqref{eq:V} with probability~$\nu$ simultaneously for all $t$ with the desired confidence level~$\kappa.$
%
%
%
%
%
% Before proving that Algorithm~\ref{alg:noise} indeed provides a sufficiently accurate noise bound, \AT{we provide a union bound by using Boole's inequality}.
% \begin{lemma}[Boole's inequality]
% \label{le:boole}
% Consider the probability space $(\tilde \Omega, \tilde{\mathcal F},\tilde \P)$ and the countable number of events $\{A_t\}_{t=1}^\infty\subset \tilde{\mathcal F}.$
%     Then, $\tilde\P[\cup_{t=1}^\infty A_i] \leq \sum_{t=1}^\infty \tilde\P[A_i].$
% \end{lemma}

We proceed by formulating the probability space in which we define our guarantees.
For Lemma~\ref{th:noise}, we already constructed a product probability space over the~$m_t$ i.i.d.\ sampled constraints.
By sampling the constraints independently over the number of iterations, we can extend the product probability space $(\Omega^{m_t}, \mathcal F^m, \P^{m_t})$ to a product probability space that is valid across all iterations~$t\geq 1$, which we refer to as the master probability space.
Due to the independence across the iterations, the master probability space can be written as $(\tilde \Omega, \tilde{\mathcal F}, \tilde\P)$, with~$\tilde\Omega \coloneqq \otimes_{t=1}^\infty \Omega^{m_t}, \tilde{\mathcal F}\coloneqq \otimes_{t=1}^\infty \mathcal F^{m_t}, \tilde\P\coloneqq \otimes_{t=1}^\infty \P^{m_t}$.
We are now ready to lift Lemma~\ref{th:noise} to hold for all iterations simultaneously.

\begin{theorem}[Simultaneous scenario bound]\label{th:simultaneous}
    Let Assumption~\ref{asm:noise_space} hold,~$\tilde\epsilon_t^{(j)}, j \in [1, m_t]$ be i.i.d.\ samples from $(\Omega, \mathcal F,\mathbb P)$, and receive~$\bar\epsilon_t$ from Algorithm~\ref{alg:noise} for all~$t\geq 1$.
    Then, for any~$\gamma,\kappa\in(0,1)$, the violation probability~\eqref{eq:V} is bounded by
$
    \tilde\P[V_t(\bar\epsilon_t)> \nu]\leq\kappa
$
\emph{simultaneously} for all~$t\geq 1$.% where~$\bar\epsilon$ is the solution of~\eqref{eq:scenario_opt} at each iteration~$t$.
\end{theorem}
\begin{proof}
The guarantees in this theorem are given with respect to the master probability measure~$\tilde\P$. % rather than~$\P^m$ to take all iterations into account.
Since~$\P^m$ can be seen as a projection of the master probability measure~$\tilde\P$, we can write the results from Lemma~\ref{th:noise} with respect to~$\tilde\P$ by fixing an iteration~$t$. 
Therefore, $\tilde\P[V_t(\bar\epsilon_t)> \nu]\leq \sum_{s=0}^{\lvert\I\rvert-1}{m_t\choose s}\nu^s(1-\nu)^{m_t-s}$ for any (fixed) $t$.
To lift this result to all iterations, we first note that the outer probability measure (in this case~$\tilde\P$) quantifies the uncertainty over time---the inner probability space governed by the measure~$\P$ does not depend on time since it is fully characterized by~$\tilde\epsilon_t^{(j)}$.
Therefore, the simultaneous bound is
%\begin{equation}\label{eq:boole}
%    \begin{aligned}
$\tilde\P[\forall t\geq 1: V_t(\bar\epsilon_t)> \nu]= \tilde\P[\cup_{t=1}^\infty V_t(\bar\epsilon_t)> \nu] \leq \sum_{t=1}^\infty\tilde\P[V_t(\bar\epsilon_t)> \nu]$,
%\end{aligned}
%\end{equation}
by Boole's inequality.
Algorithm~\ref{alg:noise} computes the number of 
~$m_t$ ($\ell.~2$) such that~$\bar\epsilon_t$ ($\ell.~4$) yields a~$\nu$-bounded violation probability~\eqref{eq:V} with the iteration-adjusted confidence level~$\kappa_t$, \ie $\tilde\P[V_t(\bar\epsilon_t)>\nu]\leq\kappa_t$, see Lemma~\ref{th:noise}.
Altogether, we have $\tilde\P[\forall t\geq 1: V_t(\bar\epsilon_t)> \nu]\leq \sum_{t=1}^\infty\kappa_t\leq\kappa$,
% Algorithm~\ref{alg:noise} requires
% \begin{align}\label{eq:kappa_t}
%     \tilde \P[V_t(\bar\epsilon_t)> \nu]\leq\kappa_t
% \end{align}
% for any~$t\geq 1$ ($\ell.~2$) with $\kappa_t\coloneqq \frac{6\kappa}{\pi^2t^2}$ (see $\ell.~1$).
% Hence,
% \begin{align}
%     \tilde\P[\forall t\geq 1: V_t(\bar\epsilon_t)> \nu] &\overset{\eqref{eq:bool}}{\leq} \sum_{t=1}^{\infty} \tilde\P[ V_t(\bar\epsilon_t)> \nu] \overset{\eqref{eq:kappa_t}}{\leq}\sum_{t=1}^\infty \kappa_t \\
%     &\leq \sum_{t=1}^\infty \frac{6\kappa}{\pi^2 t^2} = \kappa,
% \end{align}
% \AT{where the first inequality follows from Boole's inequality and the last inequality holds} since $\sum_{t=1}^\infty \frac{1}{t^2}=\frac{\pi^2}{6}$ by the solution of the Basel problem.
since~$\kappa_t\coloneqq\nicefrac{6\kappa}{\pi^2t^2}$ ($\ell.~1$) and $\sum_{t=1}^\infty \frac{1}{t^2}=\frac{\pi^2}{6}$. % by the solution of the Basel problem.
\end{proof}

Theorem~\ref{th:simultaneous} establishes that with confidence $1-\kappa$, the scenario-based bounds computed by Algorithm~\ref{alg:noise} hold uniformly across all iterations with probability $1-\nu$~\ref{con:noise}. 

\begin{comment}
\begin{remark}[On Assumption~\ref{asm:noise_space}]\label{re:noise}
\AT{
In contrast to the classic~$R$-sub-Gaussian assumption~\cite{chowdhury2017kernelized}, we assume that the noise can be sampled directly, which is not possible under the sub-Gaussian formulation since it represents a family of distributions rather than a specific one.
This assumption enables bounding the true, unknown noise with high probability (Theorem~\ref{th:simultaneous}), and arguably, the tail behavior of the noise is the most important aspect for achieving such bounds.
Since the normal distribution~$\mathcal N(0,R^2)$ attains the moment-generating-function bound of the $R$-sub-Gaussian family and therefore exhibits the heaviest admissible tails, it provides an empirical and practically convenient surrogate for implementing the most conservative case.
Assumption~\ref{asm:noise_space} thereby unifies different noise models under a single theory, allowing practitioners to incorporate prior knowledge about the noise or to employ data-driven or oracle-based representations of the aleatoric uncertainty.
}
\end{remark}
\end{comment}

% \end{remark}

\subsection{High probability confidence intervals}\label{sec:confidence}
Let us now present confidence intervals by leveraging the scenario-based noise bound derived in Theorem~\ref{th:simultaneous}~\ref{con:confidence}.

\begin{corollary}[Confidence intervals]\label{cor:confidence}
    Let the conditions in Theorem~\ref{th:simultaneous} and Assumption~\ref{asm:f} hold.
    For any~$\gamma,\kappa\in(0,1)$ and any fixed regularization constant of the GP posterior~$\eta\in(0,1]$, %simultaneously for all~$t\geq 1, a\in\domain,i\in\I$, we have
\begin{align*}
   % \tilde\P[\forall t\geq 1: \P[\lvert h_i(a)-\mu_{t,i}(a)\rvert \leq \beta_{t,i} \sigma_{t,i}(a)]\geq 1-\nu]\geq 1-\kappa,
   \tilde\P[\forall t\geq 1,i\in\I, a\in\domain: \P[h_i(a)\in C_{i,t}]\geq 1-\nu]\geq 1-\kappa,
\end{align*}
with  $C_{i,t}$ as defined in~\eqref{eq:C} and
\begin{align}\label{eq:beta}
    \beta_{t,i} \coloneqq \|h_i\|_k +\sqrt{\nicefrac{\lambda_{\max}(\Xi_t)}{\eta}} \|\bar\epsilon_{i,1:t}\|_2,
\end{align}
where~$K_t\in\mathbb R^{t\times t}$ is the covariance matrix, $\Xi_t\coloneqq K_t(K_t+\eta I_t)^{-1}$, and~$\lambda_{\max}(\Xi_t)$ denotes the largest eigenvalue of~$\Xi_t$. 
\end{corollary}
\begin{proof}
We follow the proof of~\cite[Theorem~2]{chowdhury2017kernelized}, adjusted to regularization parameters~$\eta>0$ as in \cite[Theorem~1]{Fiedler2021Practical}.
However, instead of bounding the norm of the noise vector~$\epsilon_{i,1:t}$ using a surrogate supermartingale, we use our result developed in Theorem~\ref{th:simultaneous}.
The proofs in~\cite{chowdhury2017kernelized, Fiedler2021Practical} establish that $\lvert h_i(a) - \mu_{i,t}(a) \rvert \leq (\|h_i\|_k + \sqrt{\nicefrac{1}{\eta}}\|\epsilon_{i,1:t}\|_{\Xi_t}) \sigma_t(a)$
holds deterministically, with $\|\epsilon_{i,1:t}\|_{\Xi_t}\coloneqq \epsilon_{i,1:t}^\top \Xi_t \epsilon_{i,1:t}.$
For the noise norm, we have that $\|\epsilon_{i,1:t}\|_{\Xi_t} \leq \sqrt{\lambda_{\max}(\Xi_t)} \|\epsilon_{i,1:t}\|_2$ by the Rayleigh–Ritz theorem since~$\Xi_t$ is Hermitian.
% Moreover, Theorem~\ref{th:simultaneous} implies $\tilde \P[\forall t \geq 1, \forall i\in\I: \P[\|\bar\epsilon_{i, 1:t}\|_2< \|\epsilon_{i,1:t}\|_2]> \nu]\leq\kappa$. %since~$\bar\epsilon_t$ is a non-negative vector.
%
Moreover, Theorem~\ref{th:simultaneous} ensures $\bar\epsilon_{i,t}\geq\lvert\epsilon_{i,t}\rvert$ for all~$i\in\I, t\geq 1$ with high probability.
Since~$\bar\epsilon_{i,t}$ is non-negative, this implies~$\|\bar\epsilon_{i,1:t}\|_2\geq \|\epsilon_{i,1:t}\|_2$ with high probability.
Altogether, $\tilde \P[\forall t \geq 1, \forall i\in\I: \P[\sqrt{\lambda_{\max}(\Xi_t)}\|\bar\epsilon_{i,1:t}\|_2 < \|\epsilon_{i,1:t}\|_{\Xi_t}]> \nu]\leq\kappa$, 
which implies $\tilde\P[\forall t\geq 1, i\in\I,a\in\domain: \P[h_i(a)\in [\mu_{i,t}(a)\pm\beta_{i,t}\sigma_t(a)]]\geq 1-\nu]\geq 1-\kappa$, with~$\beta_{i,t}$ as in~\eqref{eq:beta}. %, which directly proves the claim.
\end{proof}
Note that the nature of the noise enters in Theorem~\ref{th:simultaneous} when generating scenarios from the noise's probability distribution.

\section{Safe BO across noise models}
In Sec.~\ref{sec:algorithm}, we present our safe BO algorithm that accommodates observation noise under various noise models.
We prove safety and optimality of our algorithm in Sec.~\ref{sec:safety_and_optimality}.

\subsection{Algorithm}~\label{sec:algorithm}

\begin{algorithm}
\begin{algorithmic}[1]
\Require $k$, $\domain$, $S_0$, $\|h_i\|_k$, $\eta$, $\delta$, $\I$,
\For{$t=1,2,\ldots$}
\State Compute GP posterior mean~$\mu_{t,i}$ and variance~$\sigma_{t}$
\State $\bar\epsilon_t\gets$ Algorithm~\ref{alg:noise}
\State Build confidence intervals %~$C_i,t$ 
\Comment{Corollary~\ref{cor:confidence},~\eqref{eq:C}}
\If{$t>1$} $S_t\gets$~\eqref{eq:safe_set} \textbf{else} $S_t\gets S_0$
\EndIf
%\State Compute~$\ell_{i,t}, u_{i,t}, w_{i,t}\forall i\in\I$ and sets $M_t, G_t$
\State Compute~$w_{i,t}$ and sets~$M_t,G_t$
%\Comment{\eqref{eq:max}, \eqref{eq:exp}}
\State     $a_{t+1}=\argmax_{a_i\in (M_{t}\cup G_{t})}\max_{i\in\I} w_{t,i}(a_i)$% \Comment{\eqref{eq:aq}}
\If{$w_{i,t}(a_{t+1})<\delta$} break
\EndIf
\State $y_{i,t+1}\gets h_i(a_{t+1})+\epsilon_{i,t}$ \Comment{Experiment}
\EndFor
\State \Return $\argmax_{a\in S_t}\ell_{0,t}(a)$ \Comment{Optimal parameter}
\end{algorithmic}
\caption{Safe BO under general noise}
\label{alg:safe_BO}
\end{algorithm}

Algorithm~\ref{alg:safe_BO} builds up on \safeopt\ and summarizes our safe BO algorithm to solve~\eqref{eq:opt_problem}. 
First, we compute the GP mean and variance ($\ell.~2$) before receiving the scenario bound~$\bar\epsilon_t$ from Algorithm~\ref{alg:noise} ($\ell.~3$).
After building the confidence intervals ($\ell.~4$), we compute the required sets ($\ell.$~5-6)
before acquiring the next parameter~$a_{t+1}$ ($\ell.~7$).
If the uncertainty of~$a_{t+1}$ is larger than the user-defined exploration threshold~$\delta>0$ ($\ell.~8$), we conduct an experiment ($\ell.~9$) and continue the procedure.
Otherwise, we return the parameter that we estimate to safely maximize the reward ($\ell.~10$).

\subsection{Safety and optimality}\label{sec:safety_and_optimality}
Next, we prove safety and optimality of Algorithm~\ref{alg:safe_BO}.
Since exploration is restricted to the safe set~$S_t\subseteq\domain$, we can only hope to identify the maximizer of the reward~$h_0$ within a subset of the domain~$\domain$.
Thus, equivalent to~\cite{sui2015safe, berkenkamp2023bayesian}, we introduce the reachability operator $R_\delta(S)\coloneqq S \cup \cap_{i\in\Ig} \{
a\in\domain\vert \exists a^\prime\in S: h_i(a^\prime) -\delta-\|h_i\|_k d_k(a,a^\prime)\geq 0 \}.$
Moreover, we denote the reachable set by~$\bar R_\delta(S)\coloneqq \lim_{s\rightarrow\infty} R_\delta^s(S)$ with~$R^{s}_\delta(S)\coloneqq R_\delta(R_\delta^{s-1}(S))$.
Before proving that we safely converge to the~$\delta$-reachable maximizer within~$\bar R_\delta(S)$, we introduce the maximum information gain~$\gamma_{\lvert\I\rvert t}$, which bounds the information we can obtain about~$h_i$ from measurements; see~\cite{berkenkamp2023bayesian, srinivas2012information} for details.
Note that obtaining~$\gamma_{\lvert\I\rvert t}$ is non-trivial.
Nevertheless, unlike~\cite{srinivas2012information, sui2015safe, berkenkamp2023bayesian}, we do not require~$\gamma_{\lvert\I\rvert t}$ to execute Algorithm~\ref{alg:safe_BO}.
% \AT{We merely use it to derive the worst-case number of iterations for Algorithm~\ref{alg:safe_BO} to terminate and return the~$\delta$-reachable optimum.}
%Specifically, we aim for identifying the maximizer of~$h_0$ within~$\bar R_\delta (S_0)$, which is the set that results when repeatedly applying the reachability operator on the initial safe set~$S_0$.

\begin{theorem}[Safety and optimality]\label{th:safety_and_optimality}
   Let the conditions in Corollary~\ref{cor:confidence} hold, a non-empty initial safe set of parameters $S_0$ be given, and define $a^\star_t\coloneqq \argmax_{a\in S_t}\ell_{t,0}(a)$ with
\begin{align*}
    t^\star
    \geq \left\lceil\frac{8\bar\beta_{t^\star}^2\gamma_{\lvert\I\rvert t^\star} (\lvert \bar R_0(S_0) \rvert +1)}{\log(1+\eta^{-2})\delta^2}\right\rceil,
\end{align*}
%where $\lceil \zeta\rceil \coloneqq \min(s\in\mathbb N: s\geq \zeta),\forall \zeta\in\mathbb R$, 
where $\lceil \zeta\rceil$ denotes rounding~$\zeta\in\R$ up to the next integer
and~$\bar\beta_t\coloneqq \max_{i\in\I}\beta_{i,t}$.
Then, for any~$\nu,\kappa,\delta\in(0,1)$ and any~$\eta\in(0,1]$, Algorithm~\ref{alg:safe_BO} terminates after at most~$t^\star$ iterations, and the following two statements hold simultaneously with probability at least~$1-\nu$ and confidence at least~$1-\kappa$:
\begin{itemize}
    \item \emph{Safety:} $h_{i}(a_t) \geq 0, \forall i\in\Ig, \forall t\in[1,t^\star];$
    \item  \emph{$\delta$-reachable optimality:} $\ell_{0,t^\star}(a_t^\star)\geq  \max\limits_{a\in \bar{R}_\delta (S_0)}h_{0}(a)-\delta$.
\end{itemize}
\end{theorem}
\begin{proof}
%We first prove safety and then optimality.
To guarantee safety,~\cite[Lemma~11]{sui2015safe} and~\cite[Lemma~7.9]{berkenkamp2023bayesian} prove that~$\forall a\in S_t, i\in, t\geq 1, i\in \Ig: h_i(a)\geq 0$ with high probability, while~\cite[Lemma~2]{tokmak2025safe} extends this result to unknown Lipschitz constants \AT{by exploiting} Lemma~\ref{le:cont}.
We adjust~\cite[Lemma~2]{tokmak2025safe} using our confidence intervals derived in Corollary~\ref{cor:confidence} to obtain
$\tilde\P[\forall t\geq 1, \forall i\in\Ig, \forall a\in S_t: \P[h_i(a)\geq 0]\geq 1-\nu]\geq 1-\kappa$, from which safety directly follows since Algorithm~\ref{alg:safe_BO} only samples within~$S_t$.
The optimality proof proceeds analogously to the proofs of~\cite[Theorem~1]{sui2015safe} and~\cite[Theorem~4.1]{berkenkamp2023bayesian}.
The key differences from their proofs are that we rely on \AT{the RKHS norm instead of Lipschitz constant for continuity (Lemma~\ref{le:cont})} 
and that we employ the confidence intervals from Corollary~\ref{cor:confidence}, which yield guarantees holding with probability~$1-\nu$ and confidence~$1-\kappa$.
Also, we do not ensure that~$\bar\beta_t$ is non-decreasing in~$t$ by exploiting the monotonicity of the information gain.
Instead,~$\bar\beta_t$---and~$\beta_{i,t}\, \forall i \in\I$~\eqref{eq:beta}---are non-decreasing in~$t$ since $\lambda_{\max}(\Xi_t)=\lambda_{\max}(K_t(K_t+\eta I_t)^{-1})=\frac{\lambda_{\max}(K_t)}{\lambda_{\max}(K_t)+\eta}$, and~$\lambda_{\max}(K_t)$ is non-decreasing in~$t$ due to Cauchy's eigenvalue interlacing theorem.
\end{proof} 
Theorem~\ref{th:safety_and_optimality} states that Algorithm~\ref{alg:safe_BO} remains safe and terminates (see~$\ell.~8$) after at most~$t^\star$ iterations, \ie when it attains~$\delta$-reachable optimality~\ref{con:safety_optimality}.
%Note that~$t^\star<\infty$ almost surely if the samples scenarios~$\tilde\epsilon_{i,t}^{(j)}$ are finite almost surely, which is the case for standard---including unbounded---distributions.
\AT{We} not only accommodate observation noise across various noise models, but also ensure that the theoretical guarantees of Theorem~\ref{th:safety_and_optimality} hold in practice. %
%\footnote{%
%Clearly, this is only the case if the RKHS norm~$\|h_i\|_k$ (or an upper bound) is given (Assumption~\ref{asm:f}). Although commonly assumed~\cite{sui2015safe,srinivas2012information,chowdhury2017kernelized, berkenkamp2023bayesian, Bhavi2023GSO}, a priori knowledge of the RKHS norm is a delicate assumption~\cite{tokmak2024pacsbo,tokmak2025safe}.
%}
%
%
%
%
%,fiedler2024safety}.\AT{Consider removing}}
This follows from working directly with the safe set as defined in~\eqref{eq:safe_set}, while~\cite[Section~4.3]{berkenkamp2023bayesian} instead employs an alternative approach that avoids the Lipschitz constant but does so at the expense of exploration guarantees.
In fact, most BO algorithms~\cite{srinivas2012information, Bhavi2023GSO, chowdhury2017kernelized, berkenkamp2023bayesian} fix a scalar value for~$\beta_{i,t}$, which makes their safety and optimality guarantees mundane in practice. %~\cite{fiedler2024safety}.

\begin{remark}[Finite-time convergence]\label{re:beta}
\AT{For~$t^\star<\infty$, we require $\bar\beta_t\leq \mathcal O(\sqrt{t})$; a bound on the growth rate that we do not establish. % in this letter.
In fact, adversarial noise models can be chosen that ensure~$\bar\beta_t>\mathcal O(\sqrt{t})$.
We do not characterize the exact conditions under which finite-time convergence is guaranteed and only examine the growth rate of~$\bar\beta_t$ only for a specific case later in Section~\ref{sec:Franka} (Figure~\ref{fig:beta}).}
\end{remark}

\section{Experiments}
In Section~\ref{sec:toy}, we evaluate Algorithm~\ref{alg:safe_BO} using synthetic experiments with homoscedastic sub-Gaussian and heteroscedastic heavy-tailed noise.
In Section~\ref{sec:Franka}, Algorithm~\ref{alg:safe_BO} safely tunes a controller for the Franka Emika robot, while Section~\ref{sec:computational_complexity} discusses the computational complexity of the scenario-based bounds.
All experiments were conducted with hyperparameters~$\nu=10^{-1}$, $\kappa=10^{-3}$ and the Matérn32 kernel with lengthscale~$\ell=0.1$.
We compute~$\beta_{t,i}$ as in~\eqref{eq:beta} with~$\|h_i\|_k=1$.

\subsection{Synthetic examples}\label{sec:toy}
\begin{figure}
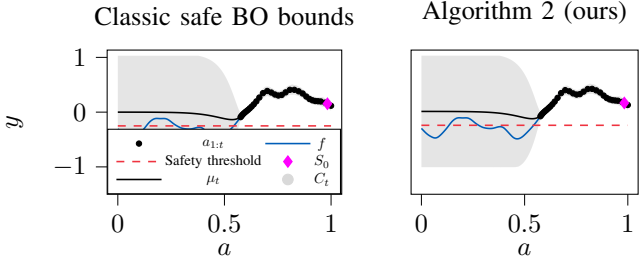

    \centering
    \input{figs/chow_uniform_final}
    \input{figs/ours_uniform_final}
    \caption{\emph{Homoscedastic sub-Gaussian noise:} Algorithm~\ref{alg:safe_BO} (ours, right) safely explores the domain to identify the maximizer with the confidence intervals specified in Corollary~\ref{cor:confidence}.
    The performance is comparable to classic safe BO algorithms (left).}
    \label{fig:Gaussian_comparison}
\end{figure}

We quantitatively compare Algorithm~\ref{alg:safe_BO} only with the standard homoscedastic sub-Gaussian safe BO bounds~\cite{chowdhury2017kernelized}, as our algorithm can operate under this standard setting while generalizing straightforwardly to heteroscedastic and heavy-tailed noise models through appropriate scenario generation.
% We first examine a toy experiment and then the Franka Emika robot.
\fakepar{Homoscedastic sub-Gaussian noise}
We
consider a scalar function~$h_i\equiv f$, \ie we have~$\lvert \I\rvert=1$.
The reward function~$f$ is a random function from the pre-RKHS of kernel~$k$, whose coefficients we scale to obtain~$\|f\|_k=1$; see~\cite[Appendix~C.1]{Fiedler2021Practical} for details.
We generate the observation noise from a homoscedastic uniform distribution with~$ \epsilon_{t} \sim \mathcal U(-10^{-3}, 10^{-3})$, scale the safety threshold to be the 40\% quantile of~$f$, and use the regularization factor~$\eta=10^{-2}$ and exploration threshold~$\delta=10^{-1}$.
Figure~\ref{fig:Gaussian_comparison} illustrates the performance of our safe BO algorithm with confidence intervals derived in Corollary~\ref{cor:confidence} (right) and standard safe BO algorithms (left).
Algorithm~\ref{alg:safe_BO} safely explores the domain and identifies the maximizer starting from the initial safe set~$S_0$. 
The performance is comparable to classic safe BO algorithms that work with sub-Gaussian confidence intervals derived in~\cite{chowdhury2017kernelized, Fiedler2021Practical} with~$R=10^{-3}$.
Next, we show that working with classic safe BO bounds outside of this noise assumption may yield safety violations.

\fakepar{Heteroscedastic heavy-tailed noise}
\begin{figure}
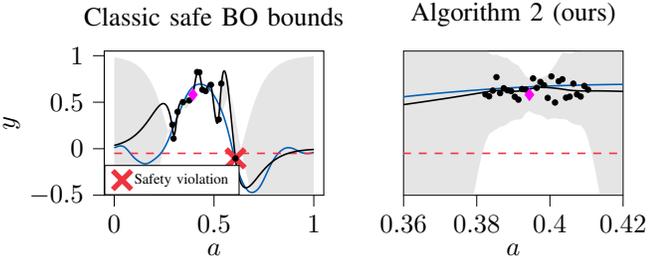

    \centering
    \input{figs/chow_heavy_final}
    \input{figs/ours_heavy_final}
    \vspace*{-0.5cm}\caption{\emph{Heteroscedastic heavy-tailed noise:} The right sub-figure provides a magnified view.
    Algorithm~\ref{alg:safe_BO} (right, ours) explores carefully and remains safe since the confidence intervals from Corollary~\ref{cor:confidence} can account for heteroscedastic heavy-tailed noise.
    In contrast, classic safe BO algorithms (left) explore more optimistically but do not account for this noise type, which yields safety violations (red cross).}
\label{fig:Student_comparison}
\end{figure}
We examine an equivalent setting as above. 
However, we here consider heteroscedastic and heavy-tailed observation noise $\epsilon_t(a)\sim\nicefrac{\lvert a\rvert}{5} \cdot \mathcal T_{10}$, where~$\mathcal T_{10}$ denotes the Student's-t distribution with ten degrees of freedom.
We further set~$R=10^{-5}$ and~$\eta=10^{-3}$.
Figure~\ref{fig:Student_comparison} compares Algorithm~\ref{alg:safe_BO} (right) with confidence intervals from Corollary~\ref{cor:confidence} to classic safe BO algorithms (left) with~$R$-sub-Gaussian confidence intervals.
By leveraging scenario-based bounds on the observation noise, Algorithm~\ref{alg:safe_BO} accommodates heteroscedastic heavy-tailed behavior, enabling more cautious exploration while maintaining safety. In contrast, classical safe BO algorithms tend to explore overly optimistically, which can ultimately lead to safety violations.

\subsection{Control parameter tuning}\label{sec:Franka}
\begin{figure}
    \centering
\includegraphics[width=7.5cm]{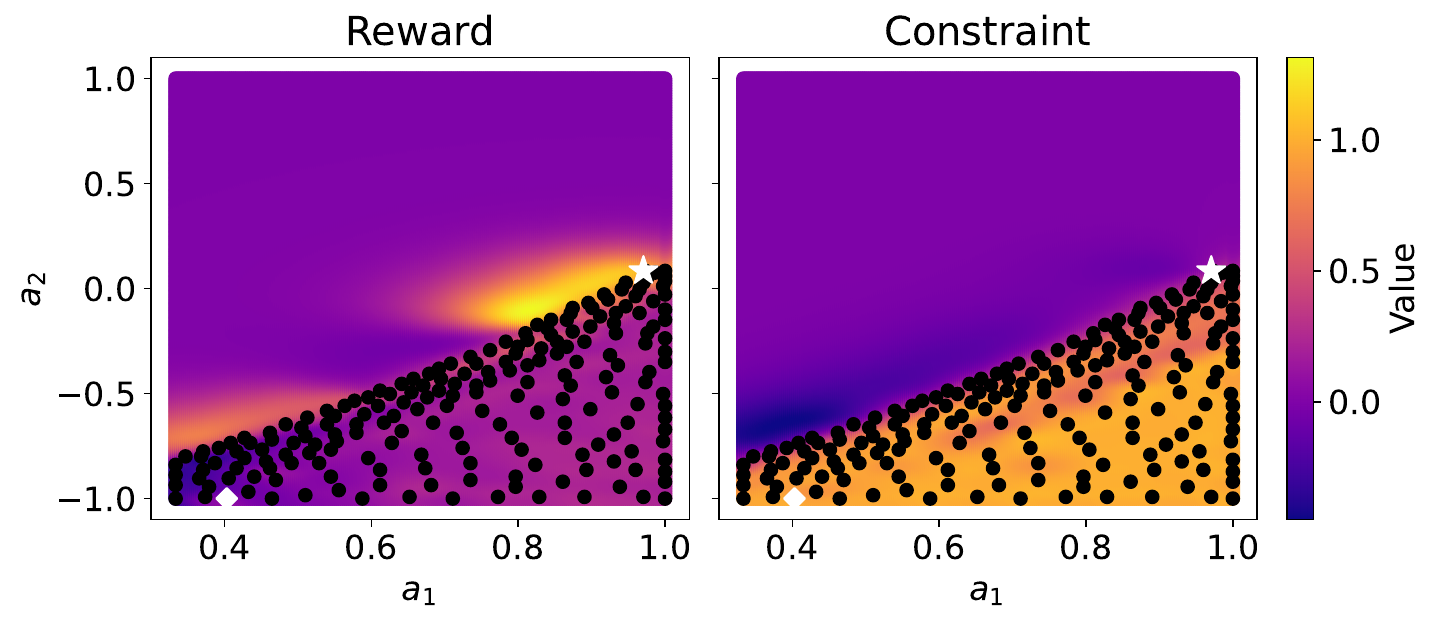}
    \caption{\emph{Franka experiment:} Exploration behavior and normalized GP mean values of the reward (left) and the constraint (right).}
    \label{fig:Franka_means}
\end{figure}

\begin{figure}
    \centering
    % This file was created with tikzplotlib v0.10.1.
\begin{tikzpicture}

\definecolor{darkgray176}{RGB}{176,176,176}

\begin{axis}[
tick align=outside,
tick pos=left,
x grid style={darkgray176},
xlabel={Iteration},
xmin=-10, xmax=210,
xtick style={color=black},
y grid style={darkgray176},
ylabel=\textcolor{aaltoBlue}{{Reward}},
ymin=-0.317924932550996, ymax=-0.185882990857901,
ytick style={color=black},
height=3cm,
width=7cm,
  ytick style={color=aaltoBlue},
  yticklabel style={color=aaltoBlue},
  y label style={color=aaltoBlue},
  y axis line style={aaltoBlue}
]
\addplot [ultra thick, aaltoBlue]
table {%
1 -0.3119230261104
10 -0.3119230261104
20 -0.274775854173919
30 -0.274775854173919
40 -0.273284618204086
50 -0.273284618204086
60 -0.273284618204086
70 -0.273284618204086
80 -0.272345506966924
90 -0.257250048610445
100 -0.248354858652466
110 -0.238520566174303
120 -0.233030801153978
130 -0.233030801153978
140 -0.215390870277704
150 -0.209303887051351
160 -0.209303887051351
170 -0.209303887051351
180 -0.191884897298497
190 -0.191884897298497
200 -0.191884897298497
};
% \addplot [thick, dashed, aaltoRed, forget plot]
% table {%
% 0 -2
% 1 -3
% };
% \addplot [semithick, aaltoRed, only marks, mark size=1, forget plot]
% table {%
% 0 10
% 1 11};
%\legend{Reward, Safety threshold, Constraint}
\end{axis}

\begin{axis}[
axis y line=right,
tick align=outside,
x grid style={darkgray176},
xmin=-10, xmax=210,
xtick pos=left,
xtick style={color=black},
y grid style={darkgray176},
ylabel=\textcolor{aaltoRed}{{Constraint}},
ymin=-0.05, ymax=1.05,
ytick pos=right,
ytick style={color=black},
yticklabel style={anchor=west},
height=3cm,
width=7cm,
  ytick style={color=aaltoRed},
  yticklabel style={color=aaltoRed},
  y label style={color=aaltoRed, yshift=2ex},
  y axis line style = {-, aaltoRed}
  ]
\addplot [semithick, aaltoRed, only marks, mark size=0.5]
table {%
0 1
1 0.971359193325043
2 0.855818510055542
3 0.942085027694702
4 0.971767127513885
5 0.848899066448212
6 0.936150312423706
7 0.986384272575378
8 0.92043924331665
9 0.858619630336761
10 0.992212116718292
11 0.962619423866272
12 0.972745776176453
13 1
14 1
15 0.881866872310638
16 0.812020659446716
17 0.991055965423584
18 0.988618552684784
19 0.983064115047455
20 0.949804604053497
21 0.997593462467194
22 1
23 0.998099029064178
24 0.989048302173615
25 0.870125412940979
26 1
27 0.987664341926575
28 1
29 0.819537699222565
30 1
31 0.844225406646729
32 1
33 1
34 1
35 1
36 1
37 1
38 1
39 1
40 1
41 1
42 1
43 1
44 0.978242576122284
45 1
46 0.930404305458069
47 0.986402332782745
48 1
49 0.645457148551941
50 1
51 1
52 1
53 0.988323867321014
54 0.586531043052673
55 1
56 1
57 0.829698503017426
58 1
59 0.985635221004486
60 0.952408254146576
61 0.991926431655884
62 0.678450465202332
63 1
64 0.987332344055176
65 0.997352838516235
66 0.922170221805573
67 0.63908725976944
68 0.985726654529572
69 0.437607288360596
70 0.534238815307617
71 0.983154833316803
72 0.990122318267822
73 0.574901342391968
74 0.648213624954224
75 0.988203048706055
76 1
77 0.969387888908386
78 0.650635063648224
79 0.849112331867218
80 0.736371576786041
81 0.717068195343018
82 0.942650437355042
83 0.847106993198395
84 0.502311050891876
85 0.856196403503418
86 0.644906759262085
87 0.685407161712646
88 0.766240835189819
89 0.84054571390152
90 0.703696489334106
91 0.60602867603302
92 0.816156327724457
93 0.304919600486755
94 0.937968015670776
95 0.759414672851562
96 0.290302217006683
97 0.776173949241638
98 0.34130647778511
99 0.661886811256409
100 0.744761943817139
101 0.518795132637024
102 0.735280692577362
103 0.670513570308685
104 0.675827145576477
105 0.741806328296661
106 1
107 0.674235463142395
108 0.602298140525818
109 0.42093500494957
110 0.719240725040436
111 0.286620914936066
112 0.694142818450928
113 0.618648886680603
114 0.727552473545074
115 0.696352064609528
116 0.241626918315887
117 0.518460988998413
118 1
119 0.684683084487915
120 0.405088394880295
121 0.947185933589935
122 0.42021632194519
123 0.120295226573944
124 0.450418949127197
125 0.188839063048363
126 1
127 1
128 0.565734624862671
129 0.246813431382179
130 1
131 1
132 1
133 0.275281012058258
134 0.39573609828949
135 0.929526269435883
136 0.861594498157501
137 0.774323046207428
138 0.383767396211624
139 0.310349464416504
140 0.880402803421021
141 0.294516026973724
142 0.902255177497864
143 0.503200173377991
144 0.305236250162125
145 0.104115135967731
146 0.331601709127426
147 0.296139001846313
148 0.360734969377518
149 0.345488250255585
150 0.976214349269867
151 0.0560372248291969
152 0.968923330307007
153 0.800584614276886
154 0.812923312187195
155 0.842805624008179
156 0.438649356365204
157 0.176006928086281
158 0.801875114440918
159 0.976176083087921
160 0.739463090896606
161 0.341069936752319
162 0.967944979667664
163 0.931734204292297
164 0.976592302322388
165 0.282266050577164
166 0.817577064037323
167 0.835395753383636
168 0.301390469074249
169 0.882772743701935
170 0.218990951776505
171 0.969270288944244
172 0.702185034751892
173 0.967400729656219
174 0.834463655948639
175 0.226400792598724
176 0.195040345191956
177 0.26625069975853
178 0.156571045517921
179 0.253637373447418
180 0.826026141643524
181 0.922495722770691
182 0.390116661787033
183 0.650217592716217
184 0.247758969664574
185 0.898970901966095
186 0.667958438396454
187 0.229270547628403
188 0.790818989276886
189 0.914152562618256
190 0.255753248929977
191 0.941628873348236
192 0.0990666002035141
193 0.644231915473938
194 0.708492398262024
195 0.63535863161087
196 0.876752078533173
197 0.194431096315384
198 0.250464409589767
199 0.603764891624451
200 0.655622124671936
};
\addplot [thick, dashed, aaltoRed]
table {%
0 0
1 0
2 0
3 0
4 0
5 0
6 0
7 0
8 0
9 0
10 0
11 0
12 0
13 0
14 0
15 0
16 0
17 0
18 0
19 0
20 0
21 0
22 0
23 0
24 0
25 0
26 0
27 0
28 0
29 0
30 0
31 0
32 0
33 0
34 0
35 0
36 0
37 0
38 0
39 0
40 0
41 0
42 0
43 0
44 0
45 0
46 0
47 0
48 0
49 0
50 0
51 0
52 0
53 0
54 0
55 0
56 0
57 0
58 0
59 0
60 0
61 0
62 0
63 0
64 0
65 0
66 0
67 0
68 0
69 0
70 0
71 0
72 0
73 0
74 0
75 0
76 0
77 0
78 0
79 0
80 0
81 0
82 0
83 0
84 0
85 0
86 0
87 0
88 0
89 0
90 0
91 0
92 0
93 0
94 0
95 0
96 0
97 0
98 0
99 0
100 0
101 0
102 0
103 0
104 0
105 0
106 0
107 0
108 0
109 0
110 0
111 0
112 0
113 0
114 0
115 0
116 0
117 0
118 0
119 0
120 0
121 0
122 0
123 0
124 0
125 0
126 0
127 0
128 0
129 0
130 0
131 0
132 0
133 0
134 0
135 0
136 0
137 0
138 0
139 0
140 0
141 0
142 0
143 0
144 0
145 0
146 0
147 0
148 0
149 0
150 0
151 0
152 0
153 0
154 0
155 0
156 0
157 0
158 0
159 0
160 0
161 0
162 0
163 0
164 0
165 0
166 0
167 0
168 0
169 0
170 0
171 0
172 0
173 0
174 0
175 0
176 0
177 0
178 0
179 0
180 0
181 0
182 0
183 0
184 0
185 0
186 0
187 0
188 0
189 0
190 0
191 0
192 0
193 0
194 0
195 0
196 0
197 0
198 0
199 0
200 0
};
\end{axis}
\end{tikzpicture}
    \caption{\emph{Franka experiment:} The maximum reward (left y-axis, blue) is increasing while the constraint value (right y-axis, red) remains above the safety threshold (dashed red line).}
    \label{fig:max_development}
\end{figure}
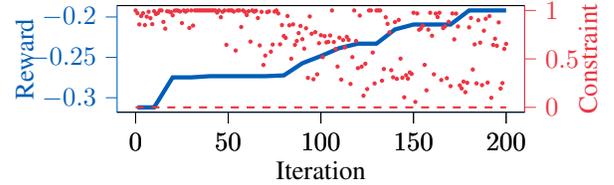

\textcolor{black}{
We next consider a set-point tracking task for the Franka Emika robot in a simulation environment based on Mujoco~\cite{todorov2012mujoco}.
Our setting is akin to~\cite[Section~5.1.1]{Bhavi2023GSO}.
We utilize the same operational space impedance controller with impedance gain $K$
and follow a linear quadratic regulator (LQR) type approach.
That is, we define the diagonal matrices~$Q$ and~$R$ as $Q=\mathrm{diag}(Q_\mathrm{r}, \kappa_\mathrm{d}Q_\mathrm{r})$, where $Q_\mathrm{r} = 10^{q_\mathrm{c}}I_3$, and $R = 10^{r-2}I_3$.  
We obtain the gain~$K$ by solving the LQR problem.
Here, we heuristically set $\kappa_\mathrm{d}=0.1$ and are then left with $q_\mathrm{c}\in [\nicefrac{1}{3}, 1]$ and $r\in [-1,1]$ as the tuning parameters.
The reward function~$f$ encourages reaching the target quickly with small inputs, whereas the constraint function~$g_1$ requires that the distance to the target decreases sufficiently.
}
%
%
%
%

%In this task, we consider the measurement noise  \AT{to} be from a homoscedastic normal distribution with~$\epsilon_{i,t}\sim\mathcal N(0,10^{-4})$.
We compute the confidence intervals (Corollary~\ref{cor:confidence}) assuming that the inherent measurement noise follows a homoscedastic normal distribution with~$\epsilon_{i,t}\sim\mathcal N(0,10^{-4})$.
We further set~$\delta=10^{-3},\eta=10^{-2}$ and start with the safe set~$S_0=\{[0.39, -1]^\top\}$.
Figure~\ref{fig:Franka_means} depicts the exploration behavior and the normalized GP means of the reward~$\mu_{0,200}$ and the constraint~$\mu_{1,200}$ after tuning the parameters for 200 iterations.
Clearly, the initial value (white diamond) yields a smaller reward than the optimal parameter after 200 iterations (white star). 
The optimization task is particularly challenging, as the region of higher rewards coincides with low constraint values.
Further, Figure~\ref{fig:max_development} illustrates that no experiment incurred a safety violation, \ie no experiment yielded a negative constraint value, and that the return of the currently believed optimum is increasing.
\begin{figure}
\centering
\begin{minipage}{0.64\columnwidth}
    % This file was created with tikzplotlib v0.10.1.
\begin{tikzpicture}

\definecolor{darkgray176}{RGB}{176,176,176}
\definecolor{steelblue31119180}{RGB}{31,119,180}

\begin{axis}[
tick align=outside,
tick pos=left,
x grid style={darkgray176},
xmin=-9.95, xmax=208.95,
xtick style={color=black},
y grid style={darkgray176},
ytick style={color=black},
width=5cm,
height=2.5cm,
xlabel={Iteration~$t$},
ylabel={$\bar\beta_t, \sqrt{t}$}
%ymode=log
]
\addplot [semithick, aaltoBlue]
table {%
0 1.00194990634918
1 1.0043112039566
2 1.00483644008636
3 1.00543117523193
4 1.0060350894928
5 1.00650334358215
6 1.00711870193481
7 1.0075751543045
8 1.00821995735168
9 1.00850176811218
10 1.00902140140533
11 1.00954210758209
12 1.00999772548676
13 1.01035368442535
14 1.01077663898468
15 1.01115155220032
16 1.01152956485748
17 1.0119765996933
18 1.01238405704498
19 1.01274371147156
20 1.01292705535889
21 1.01331126689911
22 1.01367008686066
23 1.01398146152496
24 1.01424705982208
25 1.01457977294922
26 1.01481699943542
27 1.01507532596588
28 1.01531624794006
29 1.01571714878082
30 1.01604354381561
31 1.01628446578979
32 1.01667749881744
33 1.01696991920471
34 1.01739776134491
35 1.01757705211639
36 1.01770377159119
37 1.01788210868835
38 1.01801216602325
39 1.01822853088379
40 1.01837837696075
41 1.01860773563385
42 1.01881861686707
43 1.01909387111664
44 1.01931977272034
45 1.01961255073547
46 1.01983535289764
47 1.02005982398987
48 1.02027475833893
49 1.02039194107056
50 1.02058076858521
51 1.02075433731079
52 1.02094483375549
53 1.02109169960022
54 1.02134573459625
55 1.02148151397705
56 1.0217376947403
57 1.02196311950684
58 1.02226173877716
59 1.02244639396667
60 1.02262556552887
61 1.02279186248779
62 1.02295756340027
63 1.02310395240784
64 1.02329516410828
65 1.02342700958252
66 1.02356934547424
67 1.02375876903534
68 1.0239485502243
69 1.02413523197174
70 1.02425706386566
71 1.02439177036285
72 1.02456486225128
73 1.02486538887024
74 1.02503108978271
75 1.02523446083069
76 1.02534186840057
77 1.02551448345184
78 1.02563190460205
79 1.02583241462708
80 1.02593851089478
81 1.02603304386139
82 1.02624070644379
83 1.02639830112457
84 1.02654492855072
85 1.02672946453094
86 1.02681505680084
87 1.02691102027893
88 1.02706325054169
89 1.02728772163391
90 1.02739655971527
91 1.02751684188843
92 1.02764618396759
93 1.02781438827515
94 1.02792406082153
95 1.02806484699249
96 1.02817177772522
97 1.02833700180054
98 1.02850949764252
99 1.0286431312561
100 1.02874577045441
101 1.02885496616364
102 1.02896463871002
103 1.02915835380554
104 1.02929830551147
105 1.02949810028076
106 1.02965974807739
107 1.02986919879913
108 1.03005623817444
109 1.03018939495087
110 1.03035879135132
111 1.03048622608185
112 1.03064143657684
113 1.03077697753906
114 1.03109896183014
115 1.03130686283112
116 1.03139913082123
117 1.03150629997253
118 1.0316082239151
119 1.03173971176147
120 1.03184998035431
121 1.03200435638428
122 1.03212380409241
123 1.03221881389618
124 1.03232896327972
125 1.03241920471191
126 1.03254282474518
127 1.03266727924347
128 1.03283154964447
129 1.03291738033295
130 1.03301882743835
131 1.03310716152191
132 1.03324222564697
133 1.03335225582123
134 1.03344058990479
135 1.03356850147247
136 1.0337450504303
137 1.03385376930237
138 1.03399205207825
139 1.03410542011261
140 1.03421103954315
141 1.03430378437042
142 1.0344854593277
143 1.034583568573
144 1.03468859195709
145 1.03483557701111
146 1.03497612476349
147 1.0351402759552
148 1.03540182113647
149 1.03549265861511
150 1.03557991981506
151 1.03567135334015
152 1.03578388690948
153 1.03592991828918
154 1.036017537117
155 1.03612446784973
156 1.03626596927643
157 1.0364146232605
158 1.03655350208282
159 1.03667366504669
160 1.03679263591766
161 1.0368937253952
162 1.03703629970551
163 1.03714919090271
164 1.03723120689392
165 1.03735709190369
166 1.03745710849762
167 1.03757059574127
168 1.03768396377563
169 1.03780436515808
170 1.03789699077606
171 1.03804183006287
172 1.03815186023712
173 1.03826534748077
174 1.03832995891571
175 1.03842258453369
176 1.03855431079865
177 1.03865242004395
178 1.03880250453949
179 1.03889143466949
180 1.03897142410278
181 1.03910899162292
182 1.03923571109772
183 1.03938436508179
184 1.03947007656097
185 1.03960204124451
186 1.03971827030182
187 1.03981626033783
188 1.03995025157928
189 1.04003357887268
190 1.04013562202454
191 1.04024016857147
192 1.04033851623535
193 1.04042851924896
194 1.04052305221558
195 1.04063034057617
196 1.04084253311157
197 1.04095375537872
198 1.04104471206665
199 1.0411456823349
};

\addplot[semithick, magenta, dashed, domain=0:199]
{sqrt(x)};

\end{axis}

\end{tikzpicture}
\end{minipage}\hfill
\begin{minipage}{0.31\columnwidth}
\vspace*{-0.5cm}
    \caption{\AT{Growth of confidence parameter $\bar\beta_t$ (blue) and~$\sqrt{t}$ (magenta) over iterations~$t$.}}
    \label{fig:beta}
\end{minipage}
\end{figure}
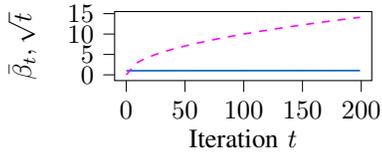

\AT{Moreover, we numerically examine the growth of~$\bar\beta_t$ in Figure~\ref{fig:beta} and observe $\bar\beta_t\leq\mathcal O(\sqrt t)$---and hence~$t^\star<\infty$ in Theorem~\ref{th:safety_and_optimality}---for this particular experiment, although a formal guarantee for this behavior in general settings is not provided in this letter; see Remark~\ref{re:beta}.}

\subsection{Computational complexity of scenario-based bounds}\label{sec:computational_complexity}

\begin{figure}
\centering
\begin{minipage}{0.61\columnwidth}
    \input{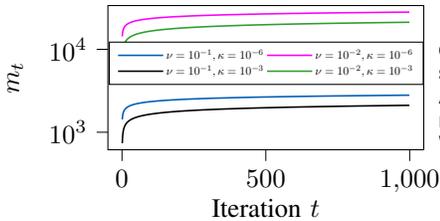}
\end{minipage}\hfill
\begin{minipage}{0.34\columnwidth}
\vspace*{-1.4cm}
    \caption{Scaling of the number of scenarios~$m_t$ in Algorithm~\ref{alg:noise} over the number of iterations~$t$. % for parameters $\{\nu=10^{-1},\kappa=10^{-6}\}$ (blue), $\{\nu=10^{-2},\kappa=10^{-6}\}$ (magenta), $\{\nu=10^{-1},\kappa=10^{-3}\}$ (black), and $\{\nu=10^{-2},\kappa=10^{-3}\}$ (green).
    We consider~$\lvert\I\rvert=1.$}
    \label{fig:scaling}
\end{minipage}
\end{figure}

\AT{Next}, we remark on the computational complexity of the scenario-based bounds we use in Algorithm~\ref{alg:noise}.
First, we examine the effects of the user-chosen accuracy parameters~$\nu,\kappa\in(0,1).$
Figure~\ref{fig:scaling} shows that the required number of scenarios~$m_t$ scales \emph{linearly} with~$\nu$.
In contrast,~$m_t$ only scales \emph{logarithmically} with~$\kappa$.
In fact, the logarithmic scaling of~$m_t$ with~$\kappa$ is a well-known result of the scenario approach~\cite{campi2018introduction}, which ensures that even extremely small confidence levels~$\kappa$ have only a mild effect on~$m_t$.
Furthermore, Figure~\ref{fig:scaling} also depicts the effect of the iteration counter~$t$ on~$m_t$.
This dependence arises due to Algorithm~\ref{alg:noise} requiring an iteration-dependent confidence level~$\kappa_t$ (Algorithm~\ref{alg:noise}, $\ell.$~1-2).
As~$t$ increases, the corresponding decrease in~$\kappa_t$ induces only a benign growth of~$m_t$,  consistent with our earlier remark on its merely logarithmic dependence on the confidence level.
Lastly, generating scenarios~$\tilde\epsilon_t^{(j)}$ in practice is usually computationally cheap if the noise is sampled from natively implemented probability distributions in standard libraries.
However, if we require importance sampling or Markov chain Monte Carlo methods, the scenario generation may become more complex.

% Finally, note that we use safe BO algorithms for control parameter tuning in an \emph{episodic} setting.
% Thus, we do not require real-time capabilities with respect to the sampling time of the closed-loop dynamics.
\AT{Finally, note that we do not require real-time capabilities with respect to the sampling time of the closed-loop dynamics since we work in an \emph{episodic} control parameter tuning setting.
}

\section{Conclusions}
We presented a safe BO algorithm reminiscent of \safeopt\ that is straightforwardly adjustable across noise models.
We proposed high-probability observation noise bounds using the scenario approach~\ref{con:noise} and integrated these bounds into frequentist confidence intervals~\ref{con:confidence}.
Using these confidence intervals, we proved safety and optimality of our safe BO algorithm~\ref{con:safety_optimality}.
Our safe BO algorithm safely explored and optimized reward functions in synthetic examples under homoscedastic sub-Gaussian and heteroscedastic heavy-tail measurement noise. 
We deployed our algorithm in tuning control parameters of safety-critical systems by tuning a controller for the Franka Emika manipulator in simulation.
\AT{Future work includes exploring data-driven noise oracles rather than analytical noise models and investigating the growth rate of~$\bar\beta_t$.}
\bibliography{IEEEabrv,references}
\bibliographystyle{IEEEtran}
\end{document}